\documentclass{amsart}
\newtheorem{theorem}{Theorem}
\theoremstyle{definition}
\newtheorem{problem}{Problem}
\newtheorem{example}{Example}

\begin{document}
\title[Problems about integral functionals]{Three topological problems 
about integral functionals on Sobolev spaces}
\author{Biagio Ricceri}
\address{Department of Mathematics, University of Catania,
Viale A. Doria 6, 95125 Catania, Italy}
\begin{abstract}
In this paper, I propose some problems, of topological nature, on the 
energy functional associated to the Dirichlet problem
$-\Delta u = f(x,u)$ in $\Omega$, $u_{|\partial \Omega}=0$.
Positive answers to these problems would produce innovative multiplicity
results on this Dirichlet problem.
\end{abstract}
\maketitle

In the present very short paper, I wish to propose some problems, of 
topological nature, on the energy functional associated to the Dirichlet 
problem
\begin{equation}\tag{$P_{f}$}
-\Delta u = f(x,u)\quad \text{in $\Omega$,}\qquad
u_{|\partial \Omega}=0.   
\end{equation}
and explain their motivations as well.

So, let $\Omega\subset \mathbb{R}^n$ ($n\geq 3$) be an open bounded 
set.
Put $X=W^{1,2}_{0}(\Omega)$.
For $q>0$, denote by $\mathcal{A}_{q}$ the class of all Carath\'eodory 
functions $f\colon \Omega\times \mathbb{R}\to \mathbb{R}$ 
such that
$$\sup_{(x,\xi) \in \Omega\times \mathbb{R}}
\frac{|f(x,\xi)|}{1+|\xi|^q} < +\infty.$$
For $0<q\leq \frac{n+2}{n-2}$ and $f\in \mathcal{A}_{q}$, put
$$\Phi_{f}(u)=
\int_{\Omega}\left ( \int_{0}^{u(x)}f(x,\xi)d\xi\right ) dx$$
and
$$J_{f}(u)=\frac{1}{2} \int_{\Omega}|\nabla u(x)|^{2}dx-\Phi_{f}(u)$$
for all $u\in X$.

So, the functional $J_{f}$ is continuously G\^ateaux differentiable on $X$ 
and one has
$$J_{f}'(u)(v) = 
\int_{\Omega}\nabla u(x)\nabla v(x)dx - \int_{\Omega}f(x,u(x))v(x)dx$$
for all $u, v\in X$. 
Hence, the critical points of $J_{f}$ in $X$ are exactly the weak 
solutions of problem $(P_{f})$.

If $q<\frac{n+2}{n-2}$, the functional $\Phi_{f}$ is sequentially 
weakly continuous, by Rellich-Kondrachov theorem. However, $\Phi_f$ may be 
discontinuous with respect to the weak topology. In this connection, 
consider the following.

\begin{example}
If $f(x,\xi)=|\xi|^{q-1}\xi$ with $0<q\leq \frac{n+2}{n-2}$, then 
$\Phi_f$ is discontinuous with respect to the weak topology.
\end{example}

In fact, if $V$ is any neighbourhood of $0$ in the weak topology of $X$, 
then $V$ does contain an infinite-dimensional linear subspace $F$ of $X$. 
Consequently, if we choose $u\in F\setminus \{0\}$, we have 
$\lambda u\in V$ for all $\lambda\in \mathbb{R}$ as well as  
$$\lim_{\lambda\to +\infty}\Phi_{f}(\lambda u) 
= \lim_{\lambda\to +\infty} 
	\frac{\int_{\Omega}|u(x)|^{q+1}dx}{q+1} \lambda^{q+1} 
= +\infty,$$
and so $\Phi_f$ is weakly discontinuous at $0$.

On the other hand, when $f$ does not depend on $\xi$, the functional 
$\Phi_f$ is weakly continuous being linear and continuous. The above 
remarks then lead to the following natural question:

\begin{problem}
Is there some $f\in \mathcal{A}_{q}$, with  $q<\frac{n+2}{n-2}$,
which is not of the form $f(x,\xi)=a(x)$, such that the functional
$\Phi_{f}$ is continuous with respect to the weak topology of $X$?
\end{problem}

To formulate the next problem, denote by $\tau_s$ the topology on $X$ 
whose members are the sequentially weakly open subsets of $X$. That is, a 
set $A\subseteq X$ belongs to $\tau_s$ if and only if for each $u\in A$ 
and each sequence $\{u_{n}\}$ in $X$ weakly convergent to $u$, one has 
$u_{n}\in A$ for all $n$ large enough.

\begin{problem}
Is there some $f\in \mathcal{A}_{q}$, with $q<\frac{n+2}{n-2}$, 
such that, for each $\lambda>0$ and $r\in \mathbb{R}$, the functional 
$J_{\lambda f}$ is unbounded below and the set $J_{\lambda f}^{-1}(r)$ 
has no isolated points with respect to the topology $\tau_{s}$?
\end{problem}

The interest for the study of Problem 2 comes essentially from the 
following result:

\begin{theorem}[{\cite[Theorem 3]{6}}]
Let $f\in \mathcal{A}_q$ with $q<\frac{n+2}{n-2}$.
Then, there exists some $\lambda^{*}>0$ such that the functional
$J_{\lambda^{*}f}$ has local minimum with respect to the topology 
$\tau_s$.
\end{theorem}

If $\Phi_f$ is weakly continuous, then the conclusion of Theorem 1
becomes stronger: the topology $\tau_s$ can be replaced by the weak 
topology. This remark is a further motivation for the study of Problem 
1.

In the light of Theorem 1, the relevance of Problem 2 is clear. Actually,
if $f$ was answering Problem 2 in the affirmative, then, by Theorem 1, for 
some $\lambda^{*}>0$, the functional $J_{\lambda^{*}f}$ would have 
infinitely many local minima in  the topology $\tau_s$. Consequently, 
problem $(P_{\lambda^{*}f})$ would have infinitely many weak solutions.

It is also worth noticing that if $f\in \mathcal{A}_{q}$ with
$q<\frac{n+2}{n-2}$ and $\lim_{\|u\|\to +\infty} J_{f}(u)=+\infty$, then 
the local minima of $J_f$ in the strong and in the weak topology of $X$ 
do coincide (\cite[Theorem 1]{3}). On the other hand, if 
$f(x,\xi)=|\xi|^{q-1}\xi$ 
with $1<q<\frac{n+2}{n-2}$, then, for some constant $\lambda>0$, it turns 
out that $0$ is a local minimum of $J_{\lambda f}$ in the strong topology 
but not in the weak one (\cite[Example 2]{3}). However, I do 
not know any example of $f$ for which $J_f$ has a local minimum in the 
strong topology but not in $\tau_s$.

To introduce the third problem (the most difficult, in my opinion), let 
me recall that in any vector space there is the strongest vector topology 
of the space (\cite[p.~42]{1}).

\begin{problem}
Denote by $\tau$ the strongest vector topology of $X$. 
Is there some $f\in \mathcal{A}_{\frac{n+2}{n-2}}$ such that the set
$\{(u,v)\in X\times X : J_{f}'(u)(v)=1\}$
is disconnected in $(X,\tau)\times (X,\tau)$?
\end{problem}

The motivation for the study of Problem 3 comes from the following 
result:

\begin{theorem}[{\cite[Theorem 1.2]{4}}]
Let $S$ be a topological space, $Y$ a real topological vector space (with
topological dual $Y^*$), and $A\colon S\to Y^{*}$ a weakly-star continuous 
operator.
Then, the following assertions are equivalent:
\begin{itemize}
\item[(i)]
The set $\{(s,y)\in S\times Y : A(s)(y)=1\}$ is disconnected.
\item[(ii)]
The set $S\setminus A^{-1}(0)$ is disconnected.
\end{itemize}
\end{theorem}

Assume that $f\in \mathcal{A}_{\frac{n+2}{n-2}}$ have the property  
required in Problem 3. Since $J_{f}\in C^{1}(X)$, clearly the operator 
$J_{f}'\colon X\to X^{*}$ is $\tau$-weakly-star continuous. Hence, by 
Theorem 2, 
the set $X\setminus (J_{f}')^{-1}(0)$ is $\tau$-disconnected. Then, this 
implies, in particular, that the set $(J_{f}')^{-1}(0)$ is not 
$\tau$-relatively compact (\cite[Proposition 3]{5}), and hence is infinite. 
So, for such an $f$, problem $(P_{f})$ would have infinitely many weak 
solutions.

Of course, to recognize the disconnectedness of the set $\{(u,v)\in 
X\times X : J_{f}'(u)(v)=1\}$ in $(X,\tau)\times (X,\tau)$, it is enough 
to check that this set is disconnected in 
$(X,\tau_{1})\times (X,\tau_{1})$, where $\tau_{1}$ is any vector 
topology on $X$ (which, to be meaningful in view of Theorem 2, should
also be stronger than the norm topology).

\providecommand{\bysame}{\leavevmode\hbox to3em{\hrulefill}\thinspace}
\providecommand{\MR}{\relax\ifhmode\unskip\space\fi MR }
\providecommand{\MRhref}[2]{%
  \href{http://www.ams.org/mathscinet-getitem?mr=#1}{#2}
}
\providecommand{\href}[2]{#2}

\end{document}